\documentclass[12pt]{article}
\usepackage{amsmath,latexsym,amssymb,amsthm,enumerate,amsmath,amscd}
\usepackage{amsthm,amsmath,amsfonts,amssymb,amscd,latexsym,stmaryrd}
\usepackage[vcentermath,enableskew]{youngtab}
\usepackage[mathscr]{eucal}
\usepackage[mathscr]{eucal}
\usepackage{bbm}
\usepackage{mathtools}
\usepackage{color}
\definecolor{med-gray}{gray}{0.5}
\definecolor{gray1}{gray}{0.85}
\definecolor{gray2}{gray}{0.90}
\definecolor{gray3}{gray}{0.64}
\definecolor{gray4}{gray}{0.55}
\definecolor{verylight-yellow}{rgb}{1,1,0.7}
\definecolor{yellow}{rgb}{1,1,0.2}
\definecolor{vivid-blue}{rgb}{0.2,0,1}
\definecolor{light-pink}{rgb}{1,0.8,1}
\definecolor{med-pink}{rgb}{1,0.6,1}
\definecolor{aqua}{rgb}{0.0, 1.0, 1.0}
\definecolor{light-gray}{rgb}{0.5, 0.9, 0.5}
\usepackage[usenames,dvipsnames]{xcolor}
\usepackage{colortbl}
\usepackage[all,cmtip]{xy}
\usepackage{comment}
\usepackage{array}
\usepackage{makeidx}
\usepackage{amssymb,multicol}
\usepackage{fancyhdr}
\usepackage{fancybox}
\usepackage{eurosym}
\theoremstyle{definition}
\usepackage{xcolor}
\usepackage{colortbl}
\usepackage[draft]{pdfpages}
\usepackage{textcomp}
\usepackage{graphicx}
\usepackage{tikz}
\usepackage{ytableau}
\usepackage{enumerate}

\theoremstyle{plain}
\newtheorem{theorem}{Theorem}[section]
\newtheorem{proposition}[theorem]{Proposition}
\newtheorem{corollary}[theorem]{Corollary}
\newtheorem{lemma}[theorem]{Lemma}
\theoremstyle{definition}
\newtheorem{definition}[theorem]{Definition}

\newtheorem{remark}[theorem]{Remark}

\newtheorem{example}[theorem]{Example}

\newtheorem{question}[theorem]{Question}
\newtheorem{algorithm}[theorem]{Algorithm}

\newtheorem{note}[theorem]{Note}

\numberwithin{equation}{section}
\numberwithin{table}{section} 
\DeclareMathOperator{\Gor}{Gor}
 
\DeclareMathOperator{\Grass}{\rm{Grass}}

\newcommand{\Hom}{\ensuremath{\mathit{Hom}}}

\newcommand{\kk}{\ensuremath{\mathsf{k}}}

\DeclareMathOperator{\Ann}{\mathrm{Ann}}

\DeclareMathOperator{\Sperner}{Sperner}

\def\rk{\mathrm{rk}}

\def\k{{\sf k}}

\def\cha{\mathrm{char}\ }

\def\G(H){\mathrm{G(H)}}
\def\Gor{\mathrm{Gor}}
\def\GGor{\mathrm{GGor}}
\def\Z(H){\mathrm{Z(H)}}
\def\ZGor{\mathrm{ZGor}}

\def\Hom{\mathrm{Hom}}

\def\ZGor{\mathrm{ZGor}}
\def\Hilb{\mathrm{Hilb}}

\def\Grass{\mathrm{Grass}}

\def\<{\left<}
\def\>{\right>}

\def\B{\mathcal{B}}

\def\ns{\footnotesize \it}

\def\Z{\mathfrak{Z}}

\def\Hess{\mathrm{Hess}}
\DeclareMathOperator{\charact}{\mathrm{char}}
\def\GGor{\mathrm{GGor}}

\def\Gor{\mathrm{Gor}}

\def\m{\mathfrak{m}}
\def\rk{\mathrm{rk}}

\newlength\distvertices

\title{Jordan type of an Artinian algebra, a survey}
\author{Nasrin Altafi\\[.05in]
{\ns Department of Mathematics and Statistics, Queen's University, Kingston, Ontario, Canada}\\[-.05in]{\ns and KTH Royal Institute of Technology, Stockholm, Sweden} 
\\[.2in] Anthony Iarrobino\\[.05in]
{\ns Department of Mathematics, Northeastern University, Boston, MA 02115,
USA.
}\\[.2in] Pedro Macias Marques\\[.05in]
{\ns Departamento de Matem\'{a}tica, Escola de Ci\^{e}ncias e Tecnologia, Centro de Investiga\c{c}\~{a}o}\\[-.05in]
{\ns  em Matem\'{a}tica e Aplica\c{c}\~{o}es, Instituto de Investiga\c{c}\~{a}o e Forma\c{c}\~{a}o Avan\c{c}ada,}\\[-.05in]
{\ns Universidade de \'{E}vora, Rua Rom\~{a}o Ramalho, 59, P--7000--671 \'{E}vora, Portugal}\\
}
\begin{document}
\date{July 3, 2023}
\maketitle
\begin{abstract} We consider Artinian algebras $A$ over a field $\k$, both graded and local algebras. The Lefschetz properties of graded Artinian algebras have been long studied, but more recently the Jordan type invariant of a pair $(\ell,A)$ where $\ell$ is an element of the maximal ideal of $A$, has been introduced.  The Jordan type gives the sizes of the Jordan blocks for multiplication by $\ell$ on $A$, and it is a finer invariant than the pair $(\ell,A)$ being strong or weak Lefschetz. The Jordan degree type for a graded Artinian algebra adds to the Jordan type the initial degree of ``strings'' in the decomposition of $A$ as a ${\sf k}[\ell]$ module.  We here give a brief survey of Jordan type for Artinian algebras, Jordan degree type for graded Artinian algebras, and related invariants for local Artinian algebras, with a focus on recent work and open problems.
\end{abstract}
\tableofcontents
\section{Introduction.}
\subsubsection{Notation.}
Let $A$ be a graded or local Artinian algebra quotient of $R=\k[x_1,\ldots,x_r]$ (polynomial ring) or of ${\mathcal R}={\sf k}\{x_1,\ldots,x_r\}$ (regular local ring) with maximal ideal $\m$ and highest socle degree $j$:  that is $A_j\not=0$, but $A_i=0$ for $i>j$.  Here, for $A$ local we take $A_i$ to be the $i$-th graded piece of the associated graded algebra $A^\ast=\bigoplus \m^i/\m^{i+1}$ of $A$.  For $A$ graded we let $\m=\oplus_1^j A_i$. The Hilbert function of $A$ is the sequence $H(A)=(h_0,h_1,\ldots, h_j)$ where $h_i=\dim_{\sf k} A_i$; the \emph {Sperner number} of $H(A)$ is the maximum value of $H(A)$. The \emph{Jordan type} $P_{\ell,A}$ of a nilpotent element $\ell\in \m$ of $A$ is the partition $P$ giving the sizes of the Jordan blocks of the (nilpotent) multiplication map $m_\ell$.   The properties of $(\ell,A)$ being strong-Lefschetz ($P=H(A)^\vee$, the conjugate of the Hilbert function viewed as a partition) or weak-Lefschetz (the number of parts of $P$ is the Sperner number) of a pair $(\ell,A)$, have been investigated as such since at least 1978 - see  \cite{St,H-W,MiNa}. Earlier, J. Brian\c{c}on in 1972 showed the strong Lefschetz property $P_{\ell,A}=H(A)^\vee$ in characteristic zero for each codimension two Artinian algebra $A$ and a generic $\ell\in R_1$ \cite{Bri}. But Jordan type is a finer concept: there are in general many partitions  that  can occur for $P_{\ell,A}$ given just the Hilbert function $H=H(A)$. A basic introductory paper is the second two authors' joint paper with C. McDaniel \cite{IMM2}; other resources include \cite{BMMN,IMM1,IMM3,CGo}.  Our attention in this note will be to the more general notion of Jordan type, as opposed to merely the Lefschetz properties.\par
Let $H$ be a sequence that occurs as the Hilbert function of an Artinian quotient of $R$ or $\mathcal R$. First, take $R$ to be the polynomial ring. We denote by $\G(H) $ and $ \GGor(H)$ the family of graded or graded Gorenstein, respectively, quotients of $R$ having Hilbert function $H$.  Now take $\mathcal R$ to be the regular local ring, and denote by $\mathrm{Z(H)}$ or $\ZGor(H)$, respectively, the family of all (not necessarily graded) quotients of $\mathcal R$ having Hilbert function $H$, or, respectively, the Gorenstein quotients of $\mathcal R$ having Hilbert function $H$. We regard these in this survey as subvarieties (not necessarily irreducible) of the Grassmanian $\Grass(R/\m^n)$, ${n=|H|}$; but some have also looked at the scheme structures, namely the Hilbert scheme $\Hilb^n(R)$  (see, for example \cite{Hui,Je,BCR,Kl,CJN} and \cite[Appendix C]{IK}). We will write $R$ for both $R$ and $\mathcal R$, when considering both at the same time. There is a natural notion of dominance of Jordan types (see Definition \ref{dominancedef}).
Our goals in this survey are 
\begin{enumerate}[(a).]
\item Review the definitions and properties of Jordan type and Jordan degree type.
\item Report on progress on the several major questions below, and
\item Suggest some further problems. 
\end{enumerate}
\subsection{Major problems.}
\label{majorsec} 
The development of the subject has been related to some questions:
\begin{enumerate}[(i).]
\item  How does Jordan type behave as one deforms the the element $\ell\in\m$, or the algebra $A=R/I$ among algebras of a given Hilbert function? Two cases: graded $A$, and local $A=\mathcal R/I$.
In particular, does the Hilbert function determine a bound (in the sense of domination) on the possible Jordan types?\par
\item For graded $A$, there is a refinement of Jordan type to a Jordan degree type \cite{IMM2}. Determine its properties and avatars
(Sections \ref{JDTsec} and \ref{avatarsec} below).\par
There is a natural generalization of Jordan type to ``contiguous Jordan type'' for graded algebras having non-unimodal Hilbert function. There are similar questions of deformation (see \cite[Section 2F, Definition 2.28ii]{IMM2}, not treated here). 
\item  Generalizations and refinements of Jordan type from graded algebras to local algebras \cite{IMS} (see Section \ref{localsec} below). 
\item When $A$ is local Gorenstein, what is the relation of these refinements of Jordan type to the symmetric decomposition of $A$ (see \cite{IM1,IMS})?
\item Using Jordan type and other invariants to show that various families  $Z(H)$ or $ZGor(H)$ have several irreducible components 
\cite{IM2} (Section \ref{3.1sec} below).  
\item Given the Artinian algebra $A$, and a fixed partition $P$ of $|A|$, what is the locus $\Z_P\subset \mathbb P(A_1)\cong \mathbb P^{r-1}$ of linear forms $\ell$ for which $P_{\ell,A}=P$? The non-Lefschetz locus \cite{BMMN}?
\item What is the relation between Jordan type and the Betti minimal resolution of $A$? \cite{Ab,AbSc}	
\item  What pairs of Jordan type partitions  $P_{\ell,A}$ and $P_{\ell^\prime,A}$ may occur together in an Artinian $A$? OR, what Jordan types $P_M$, $P_N$ may occur for a pair $(M,N)$ of $n\times n$ commuting matrices (see \cite{Kh}).
\end{enumerate}
Some of these questions are now partially answered, ideas behind them have inspired other questions that remain open. We discuss (i)-(v) in more detail below, and then pose some specific questions.\par
\subsection{What is Jordan type?}
We first present the definitions and some properties of Jordan type, and then in Section \ref{LfromJTsec} discuss the relationship to the weak and strong Lefschetz properties for graded algebras. Since the definition of Jordan type does not require grading, we start by stating it in the general setting, for a module over an algebra that may not be graded.

\begin{definition}[Jordan type]
\label{JTdef}
(See also \cite[Definition 2.1]{IMM2} and \cite[Section 3.5]{H-W})
Let $M$ be a finitely generated module over the Artinian algebra $A$, and let ${\ell\in\m}$. The Jordan type of $\ell$ in $M$ is the partition of $\dim_{{\mathsf{k}}}M$, denoted $P_\ell=P_{\ell,M}=(p_1,\ldots,p_s)$, where ${p_1\geq \cdots\geq p_s}$, whose parts $p_i$ are the block sizes in the Jordan canonical form matrix of the multiplication map ${m_\ell:M\to M}$, ${x\mapsto\ell x}$. The \emph{generic Jordan type} of $A$, denoted $P_A$, is the Jordan type $P_{\ell,A}$ for a generic element $\ell$ of $A_1$ (when $A$ is graded), or of $\m_A$ ($A$ local).
\end{definition}
The Jordan block form for the similarity class of a matrix is sometimes called the Segre characteristic, in contrast to its conjugate, the Weyr characteristic (see Note \ref{Weyrnote} below).
\begin{definition}[Jordan basis, pre-Jordan basis]
\label{Jbasisdef}
With the notation of the previous definition, a \emph{pre-Jordan basis} for $\ell$ is a basis of $M$ as a vector space over $\kk$ of the form 
\begin{equation}
\label{Jbasiseq}
{\mathcal{B}=\{\ell^iz_k \mid 1\leq k\leq s,\, 0\leq i\leq p_k-1 \}},
\end{equation}
where ${P_{\ell,M}=(p_1,\ldots,p_s)}$ is the Jordan type of $\ell$. We call the subsets $S_k=\{z_k,\ell z_k,\ldots,\ell^{p_k-1}z_k\}$ \emph{strings} of the basis $\mathcal{B}$, and each element $\ell^iz_k$ a \emph{bead} of the string. The Jordan blocks of the multiplication $m_\ell$ are determined by the strings $S_k$, and $M$ is the direct sum
\begin{equation}
\label{stringeq}
M=\langle S_1 \rangle\oplus\cdots\oplus \langle S_s \rangle.
\end{equation}
If the elements ${z_1,\ldots,z_s\in M}$ satisfy ${\ell^{p_k}z_k=0}$ for each $k$, we call $\mathcal{B}$ a \emph{Jordan basis} for $\ell$, recovering the usual definition in linear algebra, since a matrix representing the multiplication by $\ell$ with respect to $\mathcal{B}$, ordering elements as ${(\ell^{p_1-1}z_1,\ldots,z_1,\ell^{p_2-1}z_2,\ldots,z_2,\ldots,\ell^{p_s-1}z_s,\ldots,z_s)}$, is a canonical Jordan form. In that case the $\langle S_k\rangle$ are cyclic $\kk[\ell]$-submodules of $M$.
\end{definition}
The following is well-known (see \cite[\S 4.7]{Ar}, \cite{We}).
\begin{lemma}\label{preJtoJlem}
If $B$ has a pre-Jordan basis $\mathcal{B}$ as in \eqref{Jbasiseq}, then for each $k$, we have $${\ell^{p_k}z_k\in\langle \ell^az_i\mid a\ge p_k, i<k\rangle}.$$ There is a Jordan basis of $M$ derived from the pre-Jordan basis, and having the same partition invariant $P_{\ell,M}$ giving the lengths of strings.
\end{lemma}

\begin{algorithm}\label{JBalg} Often it is useful to consider a pre-Jordan basis (or a Jordan basis) to study the Jordan type of an element $\ell\in\m$. However, to compute the Jordan type of an element in a module, we do not need to choose a basis. We can consider the sequence ${(d_0,\ldots,d_{j+1})}$, where ${d_i=\dim_\kk M/\ell^iM}$, and compute the sequence of differences ${\Delta_\ell=(\delta_1,\ldots,\delta_{j+1})}$, where ${\delta_i=d_i-d_{i-1}}$. Then taking the conjugate partition of this sequence, we get the Jordan type of $\ell$ in $M$ (see \cite[Lemma 2.3]{IMM2}):
\[
P_{\ell,M}=\Delta_\ell^{\,\vee}.
\]
This is the algorithm used in Macaulay 2.
\end{algorithm}
A key  notion is specialization of Jordan types, which follows the dominance partial order on partitions
(Lemma \ref{HFdomJtlem}).
\begin{definition}[Dominance order]
\label{dominancedef} 
Let ${P=(p_1,\ldots, p_s)}$, ${p_1\ge \cdots \ge p_s}$, and ${Q=(q_1,\ldots,q_r)}$, ${q_1\ge \cdots \ge q_r}$, be two partitions of ${n=\sum p_i=\sum q_i}$. We say that $P$ dominates $Q$ (written ${P\ge Q}$, if for each ${k\in [1,\min\{s,r\}]}$, we have 
\[ 
\sum_{i=1}^k p_i\ge\sum_{i=1}^k q_i.
\]
\end{definition}
For example,  the partition $(5,4,2)\ge (5,3,2,1)$, but $(5,3,3,2)$ and $(4,4,4,1)$ are incomparable.\par
Let $P$ be a partition of $n$ we denote by $P^\vee$ the conjugate partition of $n$: switch rows and columns in the Ferrers diagram of $P$.   Let $H$ be a sequence that occurs as the Hilbert function of an Artinian algebra, and denote by $P_H$ the associated partition of $n=|A|, H^\vee$ its conjugate.
\begin{figure}[ht]
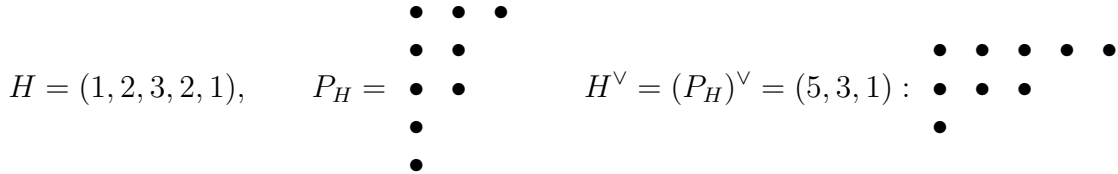

\[
H=(1,2,3,2,1),\qquad P_H=
\begin{array}{ccc}
\bullet&\bullet&\bullet\\
\bullet&\bullet&\\
\bullet&\bullet&\\
\bullet&&\\
\bullet&&
\end{array}\qquad
 H^\vee=(P_H)^\vee=(5,3,1): \begin{array}{ccccc}
\bullet&\bullet&\bullet&\bullet&\bullet\\
\bullet&\bullet&\bullet&&\\
\bullet&&&&
\end{array}
\]
\caption{Hilbert function, its partition $(3,2,2,1,1)$, and conjugate (Example \ref{1ex}).}\label{HFfig}
\end{figure}
\vskip 0.2cm

The following result is well known.
\begin{lemma}\cite[Theorem 2.5]{IMM2}\label{HFdomJtlem} Let $A$ be a standard graded Artinian algebra, and let $\ell\in A_1$ be a linear form. Then $P_{\ell,A}\le H(A)^\vee$ in the dominance partial order on partitions.
\end{lemma}
There is an analogous statement for local algebras $A$ (ibid.). \par
\begin{corollary}\label{sperncor}
 Let $A$ be an Artinian quotient of $A$ and let $\ell\in \m_A$. Then $P_{\ell,A}$ has at least as many parts as the Sperner number of $H(A)$.
\end{corollary} 
\begin{proof}That $H(A)^\vee\ge P_{\ell,A}$ and are partitions of $n=\dim_{\sf k} A$ is equivalent to $H(A)=\bigl(H(A)^\vee\bigr)^\vee\le P_{\ell,A}^\vee$ \cite[Lemma 6.3.1]{CM}. So the largest part of $H(A)$ (viewed as a partition) is less or equal the largest part in $P_{\ell,A}^\vee$, which is just the number of parts of $P_{\ell,A}$.
\end{proof}

\begin{example}[Comparison of Jordan type for algebra $B$ and associated graded algebra $A=B^\ast$]
\label{1ex}
(a) Consider the graded CI algebra $A={\k}[x,y]/I$, ${I=(x^3,y^3)=\Ann(X^2Y^2)}$, with $H(A)=(1,2,3,2,1)$ and 
$H^\vee=(5,3,1)$ (Figure \ref{HFfig}).  Here \vskip -0.2cm
\[
P_{\ell,A}=(5,3,1) \text{ for }\ell=ax+by \text{ when }ab\not=0,\text{ but }P_{x,A}=P_{y,A}=(3,3,3).
\]
The strings for ${\ell=x}$ are $\{1,x,x^2\}$, $\{y,yx,yx^2\}$, $\{y^2,y^2x,y^2x^2\}$, and $(5,3,1)>(3,3,3)$.\par\vskip 0.2cm
(b)  Consider the non-homogeneous  CI algebra $ B= \mathcal{R}/J$, with ${\mathcal{R}=\kk\{x,y}\}$, and ideal ${J=(x^3, y^3-x^2y^2)=\Ann(X^2Y^2+Y^3)}$ satisfying ${B^\ast=A}$. We have for ${\cha \k =0}$, again $P_{\ell,B}=(5,3,1)$ for $\ell=ax+by$ when $ab\not=0$, and $P_{x,B}=(3,3,3)$. But now  $P_{y,B}=(4,3,2)$, as the multiplication $m_y$ has pre-Jordan basis strings  $\{1,y,y^2,y^3=x^2y^2\}$, $\{x,xy,xy^2\}$, and  $\{x^2, x^2y\}$. Applying Algorithm~\ref{JBalg}, a Jordan basis for $m_y$ has the strings 
$\{1,y,y^2,y^3=x^2y^2\}$, $\{x,xy,xy^2\}$, and $\{x^2-y, x^2y-y^2\}$, as $y^4$, $xy^3$, and $(x^2-y)y^2$ are zero.
The algebra $B$ is a deformation of $A$, and $P_{y,B}=(4,3,2)> P_{y,A}=(3,3,3)$ in the dominance partial order, consistent with Corollary \ref{AGAcor}. 
\end{example}

The following example illustrates some of the methods of determining Jordan type for a non-homogeneous AG algebra. See also \cite[\S 2.4]{IM2}.
\begin{example}
\label{generalexample1}(Determining Jordan type, $C$ non-homogeneous.)
Let ${\mathcal{R}=\kk\{x,y,z\}}$ and ${C=\mathcal{R}/\Ann G}$, where ${G=X^3Y + Y^2Z}$. Then $C$ is a non-homogeneous AG algebra, not CI, defined by ${\Ann G=(xz,\, yz-x^3,\, z^2,\, xy^2,\ y^3)}$, with ${H(C)=(1,3,3,2,1)}$ and ${H(C)^\vee=(5,3,2)}$.\vskip 0.15cm\par\noindent
i. {\bf Generic Jordan type of $C$.} Assume ${\charact\kk\notin\{2,3\}}$ and consider a general element ${\ell\in \m_C}$. We write ${\ell=ax+by+cz+h}$, with ${h\in\m_C^{\,2}}$. Suppose ${ab\ne0}$. Then ${\ell^4=4a^3bx^3y\ne0}$. Also, ${\ell^3=a^3x^3+3a^2bx^2y+h'}$ and ${\ell^2x=a^2x^3+2abx^2y+h''}$, with ${h',h''\in\m_C^{\,4}}$ (note that ${yz=x^3}$ in $A$, so ${y^2z=x^3y\in\m_C^{\,4}}$). We can easily check that $\ell^3$ and $\ell^2x$ are linearly independent, so we have already two strings in a pre-Jordan basis for $\ell$, namely $\{1,\ell,\ell^2,\ell^3,\ell^4\}$ and $\{x,\ell x,\ell^2x\}$. According to Lemma \ref{semicontlem} the Jordan type of $\ell$ in $C$ is at most $(5,3,2)$, and we already have two string of lengths $5$ and $3$, so we will check if we can get a new string of length $2$. Note that ${\m_C^{\,3}=\langle\ell^3,\ell^2x,\ell^4\rangle}$, so if there is a further string of length two, there must be an order-one element ${\alpha\in\m_C\setminus\m_C^{\,2}}$ such that ${\ell\alpha\notin\langle\ell^2,\ell^3,\ell^4,\ell x,\ell^2x\rangle }$. Using $\ell$ and $x$ to cancel terms in $\alpha$ if necessary, we can assume that ${\alpha=z+g}$, with ${g\in\m_C^{\,2}}$. Now ${\ell\alpha=bx^3+\ell g\in\m_C^{\,3}}$, meaning there is no new length-two string. Therefore the Jordan type of $\ell$ is
\vskip -0.2cm $${P_{\ell,C}=(5,3,1,1)},$$ and since the set ${\{ax+by+cz+h\in\m_C:ab\ne0,h\in\m_C^{\,2}\}}$ is an open dense subeset of $\m_C$, this is the generic Jordan type of $C$ (Definition \ref{JTdef}). We can consider $\{z\}$ and $\{y^2\}$ as new strings to complete the pre-Jordan basis.
\vskip 0.15cm\noindent
ii. {\bf Why we cannot attain a last length-two string}.
That a last two-length string is not attainable is related to a construction from \cite[Proposition 1.33]{IM1}. The module $Q_C(1)$ can be explained by the relations between the terms $Y^2Z$ and $X^3Y$ in $G$ (we refer to \cite{IM1} for details on the $Q(a)$ modules, introduced by the second author in \cite{I2}; see also Lemma \ref{QaLemma} below). Here, $Q_C(1)$ has two homogeneous terms:

\begin{equation}
\label{generalexampleQC1}
Q_C(1)_1=\frac{(0:\m_C^{\,3})}{\m_C^{\,2}+(0:\m_C^{\,2})}
\quad\text{and}\quad
Q_C(1)_2=\frac{\m_C^{\,2}\cap(0:\m_C^{\,2})}{\m_C^{\,3}+(0:\m_C)}.
\end{equation}\noindent
Note that  
$Y^2$ is not a partial of $X^3Y$, but all further partials of $Y^2$ belong to $\langle 1,Y \rangle$, and thus are also partials of $X^3Y$. So acting on $G$ with $z$ yields ${z\circ G=Y^2}$, and this means that the class of $z$ is non-zero in $Q_C(1)_1$ (in fact, it generates this module). However, ${\m_{\mathcal{R}}z\circ G= \langle1,Y\rangle}$, so if ${\ell'\in\m_{\mathcal{R}}}$ is a lifting ot $\ell$, we have ${\ell'z\circ G=bY+d=(bx^3+dx^3y)\circ G}$, for some ${d\in\kk}$, which explains why ${\ell z\in\m_C^{\,3}}$ and its class is zero in $Q_C(1)_2$, so the module $Q_C(1)$ is acyclic. Coincidently, ${Q_C(1)=\langle z,y^2\rangle}$, so its generators are the elements we chose for the last two strings of the pre-Jordan basis.\footnote{Further examples and discussion of these points are found in \cite[\S 2.4, Remark 2.11ff.]{IM2}.}
\vskip 0.15cm\par\noindent
iii. {\bf Special Jordan types of $C$.}
When ${\ell=ax+by+cz+h}$, ${h\in\m_C^{\,2}}$ and ${ab=0}$, we find lower Jordan types in the dominance order. For instance,  
\begin{align}
\label{generalexample1lowerJT}
P_{x,C}&=(4^2,1^2), & P_{y+z,C}&=(4,2^3), & P_{y,C}&=(3^2,2^2), & P_{x^2,C}&=(2^4,1^2), &  P_{z,C}&=(2^3,1^4).
\end{align}
The strings for a pre-Jordan basis for $z$ are particularly interesting, and illustrate the issues of the non-graded case: since ${yz=x^3}$ a possible choice is $\{1,z\}$, $\{y,x^3\}$, $\{y^2,x^3y\}$, $\{x\}$, $\{x^2\}$, $\{xy\}$, $\{x^2y\}$. Note that in the strings 
$\{y,x^3\}$ and $\{y^2,x^3y\}$ there is a jump in order: the orders of $y$ and $y^2$ are $1$ and $2$, but multipliyiong by $z$ makes these orders jump to $3$ and $4$, respectively.

\smallskip\noindent
iv. {\bf Deformation $C(t)$.}
Consider the family of Artinian Gorenstein algebras $\bigl(C(t)\bigr)_{t\in\kk}$, where ${C(t)=\mathcal{R}/\Ann G(t)}$ is defined by the dual generator ${G(t)=X^3Y+Y^2Z+tYZ^2}$. Then ${C(0)=C}$, and for ${t\ne0}$, ${C(t)}$ is a CI algebra, as ${\Ann G(t)=(xz,\, ty^2-yz+x^3,\, z^2-tx^3)}$. We have ${H\bigl(C(t)\bigr)=H(C)=}$ $(1,3,3,2,1)$ for all $t$. We can check that for ${t\ne0}$ the Jordan type of ${\ell=ax+by+cz+h}$, with ${h\in\m_{C(t)}^{\,2}}$ and ${ab\ne0}$, is ${P_{\ell,C(t)}=(5,3,2)=H\bigl(C(t)\bigr)^\vee}$, admitting strings $\{1,\ell,\ell^2,\ell^3,\ell^4\}$, $\{x,\ell x,\ell^2x\}$, and $\{z,\ell z\}$. So the generic Jordan type of $C(t)$, for ${t\ne0}$, strictly dominates that of ${C=C(0)}$ which is $(5,3,1,1)$ (simply domination is required by Lemma \ref{semicontlem}). For $x$, ${y+z}$, $y$, and $x^2$, we find the same Jordan types in $C(t)$ as in \eqref{generalexample1lowerJT}, but $P_{z,C(t)}=(3^2,1^4)>P_{z,C}$.

The associated graded algebra is ${C(t)^*=R/(xz,\, ty^2-yz,\, z^2,\, x^4)}$, with ${R=\kk[x,y,z]}$, for ${t\ne0}$, and ${C(0)^*=C^*=R/(xz,\, yz,\, z^2,\, xy^2,\, y^3,\, x^4)}$. The generic Jordan type of $C(t)^*$ is the same as that of $C(t)$ as are the special Jordan types of $x$ and $x^2$, but ${P_{y+z,C(t)^*}=(3,2^3,1)}$, ${P_{y,C(t)^*}=(3,2^3,1)}$, for any ${t\in\kk}$, and ${P_{z,C(t)^*}=(2^2,1^6)}$, for ${t\ne0}$, ${P_{z,C(0)^*}=(2,1^8)}$. All these are dominated by the respective Jordan types in $C(t)$, as expected from Corollary \ref{AGAcor}.

\begin{note}\label{Weyrnote}
The \emph{Weyr characteristic}., an invariant of the similarity class of a matrix introduced by Eduard Weyr in 1885, for our nilpotent maps $m_\ell$ on $A$ is just the conjugate (switch rows and columns of the Ferrers diagram) of the Jordan partition $P_{\ell,A}$ (\cite[\S 2.4]{OCV}). See \cite{Sh} for an excellent introduction; K. O'Meara and J. Watanabe point out that for some problems the Weyr form may be more useful than the Jordan type \cite{OW}; see also \cite[Note p. 371]{IMM2} for further references.
\end{note}

\end{example}

\subsection{Historical note.} Lefschetz properties of the cohomology rings of algebraic varieties had been long studied before the algebraists adapted it.  R. Stanley showed that graded Artinian complete intersection algebras $A=R/(x_1^{a_1},\ldots, x_r^{a_r})$ satisfy a strong Lefschetz property \cite{St}: he proved this using the hard Lefschetz theorem for the cohomology of the product $\mathbb P=\mathbb P^{a_1-1}\times \cdots \times \mathbb P^{a_r-1}$ of projective spaces.  This inspired many to explore the Lefschetz properties of Artinian algebras.  Results and open problems at the time concerning Lefschetz properties of graded Artinian algebras were well set out in the 2013 foundational opus by T. Harima-J. Watanabe et al \cite{H-W} and also surveyed by J. Migliore and U. Nagel \cite{MiNa}. Other articles on the Lefschetz properties include the T. Harima articles \cite{Ha1,H} in 1995 and 1999, the 2011 B. Harbourne, H. Schenck, and A. Seceleanu on Gelfand-Tsetlin patterns and the weak Lefschetz property \cite{HSS}, and the direction of singular hypersurfaces and Lefschetz properties by R. Di Gennaro, G. Ilardi,  and J. Vall\`{e}s \cite{DIV}, a direction continued by many as E. Mezzetti, R.M. Mir\'o-Roig, G. Ottaviani on Laplace Equations and weak Lefschetz \cite{MeMO} and R.M. Mir\'o-Roig and M. Salat on Togliatti equations \cite{MR-S}. \par
Despite advocacy since 2012 at the Lefschetz conference organized by Junzo Watanabe at Tokai University, of the second author for using the finer Jordan type invariant for a pair $(\ell,A)$, it was not until \cite{IMM2} that an introduction to the topic was written. This was at the instigation of Yong-Su Shin, who asked prior to coauthoring \cite{PaSh}, where one could find an introduction to Jordan type!  There was none. The authors of \cite{IMM2} attempted to give a comprehensive introduction, including new results, doing for Jordan type what J. Migliore and U. Nagel had done earlier in the same journal in ``Tour of the strong and weak Lefschetz Properties'' \cite{MiNa}. For some topics, such as modular tensor products, they were able to exhibit several threads of work by different communities who seemed unaware of each other's work on the same subject \cite[\S 3B]{IMM2}. Several other articles by the same group treated Jordan type for certain free extensions, which are deformations of tensor products \cite[Theorem 2.1]{IMM1}; see also \cite{MCIM} which gives a connection of free extensions to invariant theory.
\par
A main advance in the study of Lefschetz properties of Artinian Gorenstein (or AG) algebras was the article of T. Maeno and J. Watanabe, showing that the ranks of multiplication by powers of a linear form $\ell$ on the degree components $A_i$ of a graded Gorenstein algebra $A$ was given by the ranks of certain higher Hessians formed from the Macaulay dual generator of $A$, at a point $p_\ell$ \cite{MW}. This result was extended and used by many, including R. Gondim \cite{Go},  Gondim and G. Zappal\`{a} \cite{GoZ1}, and it was generalized in \cite{GoZ2}, to the mixed Hessians. These have been used to prove that some Nagata idealization examples of graded AG algebras in embedding dimension at least four, are not strong Lefschetz (as \cite{CGIM}. The Hessian tools have been used recently by a growing cohort to study Jordan types for pairs $(A,\ell)$ where $A$ is a graded AG algebra and $\ell\in A_1$  (see, for example \cite{Al1,AIK,AIKY,Y1,Y2} and Section~\ref{Hesssec} below).\par

\subsubsection{Recent articles on Jordan type and Artinian algebras.}
We here mention several recent articles and research areas, with emphasis on those that mention Jordan type.
Fixing codimension two, and a Hilbert function $H$, we can study the ``Jordan cells'' $\mathbb{V}(E_P)$ of the family $G_H$, comprised of ideals having initial monomial ideal $E_P$ in a direction given by a linear form $\ell$, determined by the partition $P$, which must have ``diagonal lengths'' $T$: see \cite[Theorem 2.8]{AIKY}. The cell $\mathbb{V}(E_P)$ is comprised of all graded Artinian algebras $A=\k[x,y]/I$ such that $P_{\ell,A}=P$. The dimension of these cells, and some of their geometric properties were known \cite{Y1,Y2,IY}; the article \cite{AIKY} determines the
generic number of generators of ideals in each cell \cite[Theorems 3.11, 5.15]{AIKY} using a decomposition of cells into a product of simpler components.  See Question \ref{CellsZHquest}.
\par
There has been the beginning of tying the Jordan type with the Betti resolution of $A$, see N. Abdallah and H. Schenck \cite{AbSc} and N. Abdallah is \cite{Ab}, and as well J. Jelisiejew, S. Masuti and M. Rossi's  \cite{JMR}, where they investigate local complete intersections of codimension three, also the book-length \cite{KKRSSY} has some Betti vs. Jordan type calculations. \par
Not totally unrelated, the preprint  \cite{AAIY} studies Jordan types for codimension three graded Gorenstein algebras of Sperner number at most $5$ and all linear forms. This is facilitated by the D. Buchsbaum-D. Eisenbud Pfaffian structure theorem and related work \cite{BuEi,CV,Di} which specifies the Betti resolutions possible given $H(A)$. The results  are still complex with $26$ Jordan types for $H=(1,3,4^k,3,1)$ when $k\ge 3$ and $47$ for $H=(1,3,5^4,3,1)$.\par
 In \cite{A-N} the weak Lefschetz property and Jordan types for linear forms of a class of graded AG algebras, called Perazzo algebras, of codimension five were studied. For Perazzo algebras, the multiplication map $\ell^{j-2}$ from degree $1$ to degree $j-1$ does not have maximal rank, where $j$ is the socle degree. Thus, the strong Lefschetz property for this family is not satisfied. In \cite{A-N} all Jordan types for linear forms of Perazzo algebras of codimension five with  the smallest possible Hilbert function were determined.
\section{Properties of Jordan type, and of Jordan degree type.}

\subsection{Lefschetz properties and Jordan type.}
\label{LfromJTsec}

\begin{definition}[Lefschetz properties]\label{WLdef}
Let $A$ be a graded Artinian algebra of highest socle degree $j$ and let ${\ell\in A_1}$. We say that the pair $(A,\ell)$ is a \emph{weak Lefschetz} (WL) if for each ${i\ge0}$ the multiplication map ${\times\ell:A_i\to A_{i+1}}$ has maximal rank. The algebra $A$ satisfies the \emph{weak Lefschetz property} (WLP) if it has a WL element. We say that the pair $(\ell,A)$ is \emph{strong Lefschetz} (SL) if for each ${i,d\ge0}$ the multiplication map ${\times\ell^d:A_i\to A_{i+d}}$ has maximal rank. The algebra $A$ satisfies the \emph{strong Lefschetz property} (SLP) if it has a SL element.
\end{definition}
The following result part A is a portion of \cite[Prop. 2.10]{IMM2}; part B is essentially \cite[Lemma 2.11]{IMM2}, shown when $H(A)$ is also symmetric in \cite[Prop. 3.5]{H-W}. We say that a Hilbert function $H(A)=(h_0,h_1,\ldots, h_j)$ is 
\emph{unimodal} if there is an integer $k$ such that $h_0\le \cdots \le h_k$ and $h_k\ge h_{k+1}\ge \cdots \ge h_j$. Recall that the Sperner number Sperner$(A)$ is the maximum value of $H(A)$.
\begin{lemma}\label{SLWLlem}A. Let $A$ be a graded Artinian algebra (possibly non-standard), and $\ell\in A_1$. Then the following are equivalent
\begin{enumerate}[(i.)]
\item The pair $(A,\ell)$ is strong Lefschetz;
\item  The Jordan type $P_{A,\ell}=H(A)^\vee$, the conjugate of the Hilbert function viewed as a partition.
\end{enumerate}
B.  Assume further that $H(A)$ is unimodal. Then the following are equivalent
\begin{enumerate}[(i).]
\item  The pair $(A,\ell)$ is weak Lefschetz.
\item  The dimension $\dim_k A/\ell A= \Sperner(A)$.
\item  The number of parts of the Jordan partition $P_{A,\ell}$ is $\Sperner(A)$, the minimum possible given $H(A)$ (Corollary \ref{sperncor}).
\end{enumerate}
\end{lemma}
\begin{proof} The proof of Lemma \ref{SLWLlem}(A) under the hypothesis is a bit subtle see \cite[Prop. 2.10]{IMM2}.  For Lemma \ref{SLWLlem}(B), the proof of
${\mathrm{B(i)}\Leftrightarrow \mathrm{B(ii)}}$ is straightforward from the definitions;  the proof of  
${\mathrm{B(ii)}\Leftrightarrow \mathrm{B(iii)}}$ follows from decomposing $A$ as a direct sum of \emph{strings}
(Lemma \ref{preJtoJlem}).\par
\end{proof}\par

\subsection{Higher Hessians and mixed Hessians.}\label{Hesssec} 

Graded Artinian Gorenstein algebras are determined by a single polynomial in the Macaulay dual ring, by a result of F.H.S. Macualay \cite{Mac2}. Let $A=R/\Ann F $ be an Artinian Gorenstein algebra with dual generator $F\in \mathcal{E}_j=\kk[X_1,\ldots, X_r]_j$, where $\Ann F$ is the ideal generated by all the forms $g\in R$ such that $g\circ F=0$. T. Maeno and J. Watanabe \cite{MW} introduced higher Hessians associated to the dual generator $F$ and provided a criterion for Artinian Gorenstein algebras having the SLP. \par We first briefly recall the Macaulay duality \cite{Mac3}, see \cite[\S 21.2]{Ei}, \cite{I3}; the recent emendation by J. O. Kleppe and S. Kleiman gives a geometric view consistent with studying deformation \cite{KK}. We let $R=\kk[x_1,\ldots,x_r]$ act on $\mathcal{E}$ by contraction\footnote{When $\cha \k=0$ or $\cha \k>j$ we may use the usual differentiation action, see \cite[Appendic A]{IK}.} where for ${u\ge k}$, ${x_i^k\circ X_j^u=\delta_{i,j}X_i^{u-k}}$ (we will call this ${\partial^k/\partial X_i^k \circ X_j^u}$) and extending this multilinearly to an action of $h\in R$ on $F\in\mathcal E$. 
\begin{equation} 
h\circ F=h(\partial/\partial X_1,\ldots, \partial/\partial X_r)\circ F,
\end{equation}
so taking $F=X_1^{\,3}X_2^{\,2}+X_1X_2^{\,4}$ we have $x_1x_2^{\,2}\circ  F=X_1^{\,2}+X_2^{\,2}$.
\begin{definition}\cite[Definition 3.1]{MW}
Let $F$ be a polynomial in $\mathcal E$ and $A= R/\Ann F$ be its associated Gorenstein algebra. Let $\mathcal{B}_{k} = \lbrace \alpha^{(k)}_i\rbrace_i$ be an ordered $\kk$-basis of $A_k$. The entries of the $k$-th Hessian matrix of $F$ with respect to $\B_k$ are 
\[
\Hess^k(F)=(\alpha^{(k)}_u\alpha^{(k)}_v F)_{u,v}.
\]
\end{definition}
\noindent 
Note that when $k=1$, $\Hess^1(F)$ coincides with the usual Hessian. P. Gordan and M. Noether proved that the (first) Hessian of every homogeneous form $F$ in at most $4$ variables has non-zero determinant unless $F$ defines a cone \cite{GorNo}. This is no longer the case in polynomial rings with at least $5$ variables: a family of forms that do not define a cone and for which the Hessian has zero determinant was provided by \cite{GorNo} and \cite{Per}, they are called Perazzo forms.  \par \noindent 
Up to non-zero constant multiple, $\det \Hess^k(F)$ is independent of the choice of basis $\B_k$.  For every $0\le k\le \lfloor\frac{j}{2}\rfloor$ and a linear form $\ell=a_1x_1+\cdots +a_rx_r$ the rank of $\times \ell^{j-2k}:A_k\rightarrow A_{j-k}$ is equal to the rank of $\Hess^k_\ell(F)$; i.e.\ the Hessian matrix evaluated at the point $P_\ell= (a_1,\dots ,a_r)$ -- see Theorem \ref{rankthm} below.  For now we state, 
\begin{theorem}\cite[Theorem 3.1]{MW},\cite{Wat2}.
An AG algebra $A = R/\Ann F$ with socle degree $j$ has the SLP if and only if there exists linear form $\ell\in R_1$ such that 
$$
\det\Hess^k_\ell(F)\neq 0,
$$
for every $k=0,\dots , \lfloor\frac{j}{2}\rfloor$.
\end{theorem}

As mentioned above, for Perazzo forms $F$ the determinant of the first order Hessian, $\Hess^1(F)$, is identically zero. So by the above theorem the associated AG algebra of a Perazzo form fails to have the SLP. The WLP and Jordan types of  Perazzo algebras in $5$ variables have been studied in \cite{FMM} and \cite{A-N}.

For an AG algebra for which all higher Hessians have non-vanishing determinants, the above theorem shows that for a general enough linear form $\ell$  all the multiplication maps $\ell^{j-2k}: A_k\to A_{j-k}$ have maximal rank. It is natural to ask: if an AG algebra A has at least one Hessian with vanishing determinant, which multiplication maps have maximal rank and which ones do not?
R. Gondim and G. Zappal\`{a} \cite{GoZ2} introduced mixed Hessians that generalize the higher Hessians.

\begin{definition}\cite[Definition 2.1]{GoZ2}
Let $A= R/\Ann F$ be the AG algebra associated to $F\in \mathcal{E}_j$. Let $\mathcal{B}_{k} = \lbrace \alpha^{(k)}_i\rbrace_i$ and $\mathcal{B}_{u} = \lbrace \beta^{(u)}_i\rbrace_i$ be $\kk$-basis of $A_k$ and $A_u$ respectively. The entries of the \emph{mixed Hessian} matrix of order $(k,u)$ for $F$ with respect to $\mathcal{B}_k$  and $\mathcal{B}_u$ is given by
$$
\Hess^{(k,u)}(F)=(\alpha^{(k)}_u\beta^{(u)}_v F)_{u,v}.
$$
\end{definition}
\noindent Notice that this generalizes the definition of higher Hessians and we have $\Hess^k(F)=\Hess^{(k,k)}(F)$.
\begin{theorem}\label{rankthm}\cite[Theorem 2.4]{GoZ2}. Let $A$ be an AG standard graded algebra. Then the rank of the mixed Hessian matrix of order $(k,u)$ evaluated at the point $P_\ell= (a_1,\dots ,a_r)$ is the same as the rank of the multiplication map $\ell^{u-k}:A_k\rightarrow A_u$ for $\ell=a_1x_1+\cdots +a_rx_r$.
\end{theorem}

The method of higher Hessians and mixed Hessians has been used to study the Lefschetz properties for graded AG algebras in the literature, for instance see \cite{Al1,Go,FMM}. The ranks of higher and mixed Hessians together at the point $P_\ell$ completely determine the ranks of multiplication maps by different powers of the linear form $\ell$ in all degrees, and hence, when the Hilbert function $H(A)$ is unimodal, the Jordan degree type of $(\ell,A)$ (Proposition \ref{JDTrankprop}). B. Costa and R. Gondim  in \cite{CGo} determined the Jordan types for general linear forms of AG algebras having low codimension and low socle degree in terms of the ranks of the associated mixed Hessians.

The first and second authors with L. Khatami classified \cite{AIK} all partitions that can occur as Jordan types of linear forms for AG algebras in codimension two (these are exactly complete intersection algebras by \cite{Mac1}) having a fixed Hilbert function.
It has been shown that in codimension two, the Jordan types of linear forms of AG algebras are completely determined by the rank of higher Hessians. In fact, they are uniquely determined by the sets of higher Hessians that have vanishing determinants.
\begin{theorem}\cite[Theorem 3.8]{AIK}
Assume that $H=\left( 1,2,3, \dots , d^k,\dots ,3,2,1\right)$, is a Hilbert function of some complete intersection algebra for $d\geq 2$ and $k\geq 2$ $(k=1$, respectively). Let $P$ be a partition that can occur as the Jordan type of a linear form and an Artinian complete intersection algebra having Hilbert function $H$. Then the following are equivalent.
\begin{enumerate}[(i).]
\item $P=P_{\ell,A}$ for a linear form $\ell\in R_1$ and an Artinian complete intersection algebra $A=R/\Ann F$, and there is an ordered partition $n=n_1+\cdots +n_c$ of an integer $n$ satisfying $0\leq n\leq d$ (or $0\leq n\leq d-1$, respectively) such that $\det \Hess^{n_1+\cdots +n_i-1}_\ell(F)\neq 0$, for each $1\le i\le c$, and the remaining Hessians are zero;
\item $P$ satisfies
\begin{equation}\label{JTpartition}
P=\big(p_1^{n_1}, \ldots, p_c^{n_c}, (d-n)^{d-n+k-1}\big),
\end{equation}
 where 
$p_i=k-1+2d-n_i-2\displaystyle(n_1+\cdots +n_{i-1})$, for $1\leq i \leq c$.
\end{enumerate}
\end{theorem}
The above theorem shows that given $H$ there are exactly $2^d$ (when ${k>1}$), or $2^{d-1}$ (when ${k=1}$) possible Jordan types for Artinian complete intersection algebras of codimension two with Hilbert function $H$. These correspond to the vanishing subsets of higher Hessians, $\Hess^i(F)$ for $i=0, \ldots, d-1$. 

In codimension two, the ranks of mixed Hessians $\Hess^{(k,u)}(F)$ for every ${k\neq u}$ are determined by the ranks of higher Hessians $\Hess^{k}(F)$ for $0\le k\le d-1$ \cite[Proposition 3.12]{AIK}. This is no longer true in higher codimensions. The first author in \cite{Al2} introduced rank matrices associated to AG algebras and linear forms which provide same information as the ranks of mixed Hessians. The first and second authors  with N. Abdallah and J. Yame\'ogo in a work in progress study all partitions that can occur as Jordan types for AG algebras of codimension three and low Sperner numbers (maximum value of the Hilbert function) \cite{AAIY}.  

\subsection{Jordan degree type.}\label{JDTsec}

If $A$ is a graded algebra, quotient of the polynomial ring $R$, and $M$ is a graded module over $A$ (possibly non-standard graded), we can consider a refinement of Jordan type, called \emph{Jordan degree type} (see \cite[Definition 2.28]{IMM2}, where several equivalent ways of defining it are presented and different notations are proposed).

\begin{definition}
Let $A$ be a graded algebra and let $M$ be a finite graded module over $A$. Let ${\ell\in A_k}$, ${k\ge 1}$ be a homogeneous element. If $\mathcal{B}$ is a Jordan basis for the multiplication by $\ell$ in $M$ as in \eqref{Jbasiseq} and all the elements of $\mathcal{B}$ are homogeneous then the \emph{Jordan degree type} of $\ell$ in $M$ is the sequence of pairs of integers
\[
\mathcal{S}_{\ell,M}=\bigl( (p_1,\nu_1),\ldots,(p_s,\nu_s) \bigr),
\]
where $\nu_k$ is the degree of the first bead $z_k$ in the string $S_k=\{z_k,\ell z_k,\ldots,\ell^{p_k-1}z_k\}$, and, by reordering the strings if necessary, we require that if ${p_k=p_{k+1}}$ for some $k$ then ${\nu_k\le\nu_{k+1}}$.\par
We may write the Jordan degree type in list manner:
\begin{equation}\label{listJDTeq}
\mathcal S_{\ell,A}=\sum_{p}\sum_\nu (p,\nu)^{\eta(p,\nu)}, 
\end{equation}
where $\eta(p,\nu)$ is the multiplicity of $(p,\nu)$, and the sum is over distinct pairs $(p,\nu)$.\par
For convenience we may write the pair $(p,\nu)$ as $p_{\nu}$, and abbreviate a list of pairs having consecutive degrees $(p,\nu),(p,\nu+1),\ldots,(p,\nu+e)$ as $(p\uparrow_{\nu}^{\nu+e})$ in listing Jordan degree type -- e.g.\ the list $(5,2),(5,3),(5,4)$ could be written as $(5\uparrow_{2}^{4})$.\par
\end{definition}

\begin{example}
\label{generalexample2}
Consider the family of graded Artinian algebras $\bigl(C(t)^*\bigr)_{t\in\kk}$ associated to $C(t)$ from Example~\ref{generalexample1}iv. For ${t\ne0}$, the Jordan degree type of any linear form ${\ell=ax+by+cz}$, with ${ab\ne0}$, satisfies ${\mathcal{S}_{\ell,C^*}=\bigl( (5,0),(3,1),(2,1)\bigr)},$
which is the generic JDT of $C(t)^*$, and is the JDT associated with the Hilbert function ${H=(1,3,3,2,1)}$ (see Note \ref{HFJDTnote}). If ${t=0}$,  ${C^*=C(0)}^*$ has generic JDT ${\mathcal{S}_{\ell,C^*}=\bigl( (5,0),(3,1),(1\uparrow_{1}^{2})\bigr)}$.
A few examples of special Jordan degree types are ${\mathcal{S}_{x,C(t)^*}=\bigl( (4\uparrow_{0}^{1}),(1\uparrow_{1}^{2})\bigr)}$ and ${\mathcal{S}_{y,C(t)^*}=\bigl( (3,0),(2\uparrow_{1}^{3}),(1,1)\bigr)}$, for any ${t\in \kk}$, and ${\mathcal{S}_{z,C(t)^*}=\bigl( (2\uparrow_{0}^{1}),(1\uparrow_{1}^{4}),(1\uparrow_{2}^{3})\bigr)}$, for ${t\ne0}$, while 
${\mathcal{S}_{z,C(0)^*}=\bigl( (2,1),(1\uparrow_{1}^{3}),(1\uparrow_{2}^{4}),(1\uparrow_{1}^{2})\bigr)}$.
\end{example}

Jordan degree type also enjoys a semi-continuity property. To be able to state it, we need to define a Hilbert function associated with a Jordan degree type and a partial order among Jordan degree types sharing the same Hilbert function. We will adopt here purely combinatorial definitions, without requiring that a given sequence is the Jordan degree type of some linear form in a module.

\begin{definition}
Let ${\mathcal{S}=\bigl( (p_1,\nu_1),\ldots,(p_s,\nu_s) \bigr)}$ be a sequence of pairs of non-negative integers satisfying 
\begin{equation}
\label{asJDTsequence}
p_1\ge\cdots\ge p_s\qquad\text{and}\qquad {\nu_k\le\nu_{k+1}} \text{ whenever } {p_k=p_{k+1}}.
\end{equation}
Then 
\begin{enumerate}[(i.)]
\item The \emph{partition associated to} $\mathcal{S}$ is ${P=(p_1,\ldots,p_s)}$. 
\item The \emph{Hilbert function associated to} $\mathcal{S}$ is that defined by
\[
H(\mathcal{S})_i=\#\{k:\nu_k\le i<\nu_k+p_k\}.
\]
\item For each ${i\ge0}$, we define the \emph{truncation} $\mathcal{S}_{\le i}$ as the sequence we obtain from $\mathcal{S}$ by omitting the pairs for which ${\nu_k>i}$, then replacing each other pair $(p_k,\nu_k)$ in $\mathcal{S}$ by the pair $(\min\{p_k,i+1-\nu_k\},\nu_k)$.
\item We say that two pairs $(p_k,\nu_k)$, $(p_l,\nu_l)$, can be \emph{combined} or \emph{concatenated} if ${\nu_l=\nu_k+p_k}$. In this case, the result of the concatenation is the pair $(p_k+p_l, \nu_k)$. 
\end{enumerate}
\end{definition}
\begin{remark}
Note that if $\mathcal{S}$ is the Jordan degree type of a linear form ${\ell}$ in a finite graded module $M$ over an Artinian algebra $A$ then ${H(\mathcal{S})=H(M)}$. Moreover, the truncation $\mathcal{S}_{\le i}$ is the Jordan type of $\ell$ in the module ${M/\m^{i+1}M}$.
\end{remark}
The following definition is adapted from \cite[Definition 2.28 (iii) and (vi)]{IMM2}.
\begin{definition}[Concatenation partial order and dominance partial order]\label{podef}
Let $\mathcal{S}$ and $\mathcal{S}'$ be two
sequences of pairs of non-negative integers satisfying \eqref{asJDTsequence}, with ${H(\mathcal{S})=H(\mathcal{S}')}$. 
\begin{enumerate}[(i.)]
\item Concatenation partial order: we say that ${\mathcal{S}\ge_c\mathcal{S}'}$ if $\mathcal{S}$ can be obtained from $\mathcal{S}'$ by a sequence of concatenations. 
\item Dominance order for JDT: we say that ${\mathcal{S}\ge\mathcal{S}'}$ if for each ${i\ge0}$, the partition associated to $\mathcal{S}_{\le i}$ is greater or equal to that associated to $\mathcal{S}'_{\le i}$, in the dominance partial order (Definition \ref{dominancedef}). 
\end{enumerate}
\end{definition}
\noindent
In \cite[Definition 2.28 (vi)]{IMM2}, the dominance partial order is defined only for Jordan degree types sharing the same partition. Here we have adopted an extension of this order, comparing sequences which may have different associated partitions.
\begin{remark}
Note that if ${\mathcal{S}\ge_c\mathcal{S}'}$ then ${\mathcal{S}\ge\mathcal{S}'}$. To see this, suppose that $\mathcal{S}$ is obtained from $\mathcal{S}'$ by one concatenation, of the pairs $(p_k,\nu_k)$ and $(p_l,\nu_l)$. If ${i<\nu_l}$ then  ${\mathcal{S}_{\le i}=\mathcal{S}'_{\le i}}$. If ${i\ge\nu_l}$ then the pair $(\min\{p_k+p_l,i+1-\nu_k\},\nu_k)$ occurs in ${\mathcal{S}_{\le i}}$, while the pairs $(p_k,\nu_k)$ and $(\min\{p_l,i+1-\nu_l\},\nu_l)$ occur in ${\mathcal{S}'_{\le i}}$, and this is the only difference between ${\mathcal{S}_{\le i}}$ and ${\mathcal{S}'_{\le i}}$. In any case, we see that the partition associated to $\mathcal{S}_{\le i}$ is greater or equal to the one associated to $\mathcal{S}'_{\le i}$. The result follows by transitivity.
\end{remark}
\begin{example}
Let ${\mathcal{S}=\bigl( (4,0),(2,1) \bigr)}$ and  ${\mathcal{S}=\bigl( (3,0),(3,1) \bigr)}$. Then ${\mathcal{S}\ge\mathcal{S}'}$, but $\mathcal{S}$ and $\mathcal{S}'$ are not comparable under the concatenation partial order. This shows that the dominance partial order is a proper extension of the  concatenation partial order.
\end{example}

\begin{lemma}
\cite[Lemma 2.29]{IMM2}
\label{JDTdeformlem}
Fix  
${\ell\in R_1}$. Let $\bigl(A(w)\bigr)_{w\in W}$ 
be a family of graded Artinian algebras with constant Hilbert function 
${H\bigl(A(w)\bigr)=H}$, and constant Jordan type ${P_{\ell, A(w)}=P_\ell}$. Let ${w_0\in W}$ be such that for ${w\in W\setminus\{w_0\}}$, also the Jordan degree type is constant, ${\mathcal{S}_{\ell,A(w)}=\mathcal{S}_{\ell}}$. Then ${\mathcal{S}_{\ell}\ge\mathcal{S}_{\ell,A(w_0)}}$ in the dominance partial order of Definition \ref{podef}(ii).
\end{lemma}
\begin{note}\label{HFJDTnote}  We may assign a contiguous Jordan type and JDT$(H)$ to a Hilbert function: we just regard $H$ as a bar graph, and list the rows of the bar graph with their initial degrees \cite[Definition 2.8ii]{IMM2}. In \cite[Prop. 2.32]{IMM2} is shown a concatenation  inequality JDT$(H)\ge_c \mathcal S_{\ell,A}$ for a graded $A$ of Hilbert function $H$. \par However, using our extended
Definition \ref{podef}ii. above of dominance of JDT, we have also that for $A$ graded with $H(A)=H$,
\begin{equation}
\text{JDT }(H)\ge \mathcal S_{\ell,A},
\end{equation}
in the dominance partial order.  This follows from the inequality $(H_{\le i})^\vee\ge P_{\ell,A/\m_A^{i+1}}$, the Jordan type of $A/\m_A^{i+1}$, for each $i\in [0,j]$, and the equivalence between JDT and the sequential Jordan type  (Lemma \ref{SJT=JDTlem}).
\end{note}

\begin{remark}
\label{JDTvsJTrem}
Note that Lemma \ref{JDTdeformlem} implies an analogous result for the Jordan type in a family of graded algebras of fixed Hilbert function.  This is generalized in Proposition \ref{SJTdeformprop} to families of possibly non-graded algebras.\par
The Jordan degree type for a graded Artinian algebra is a finer invariant than the \linebreak Jordan type \cite[Example 3.1]{IMM2}. Taking  $B=\mathsf{k}[x,y,z]/(yz,x^2y,xy^2,z^3,x^4,y^4)$ and  $A=\mathsf{k}[x,y,z]/(y^2,x^2z,x^2y, z^3,x^6)$, we have $H(A)=H(B)=(1,3,5,4,2,1)$ and  $P_{z,A}=P_{z,B}=(3^4,1^4)$; but $\mathcal{S}_{z,A}=\bigl((3\uparrow_0^2),3_1,(1\uparrow_2^5)\bigr)$ and $\mathcal{S}_{z,B}=\bigl((3\uparrow_0^3),(1\uparrow_1^3),1_2\bigr)$.

The first author in \cite[Example 4.2]{Al3} provided AG algebras $A$ and $B$, quotients of $R$, having Hilbert function $(1,3,6,9,9,9,6,3,1)$ and a linear form $\ell\in R_1$
for which the Jordan types of pairs $(A,\ell)$ and $(B,\ell)$ are equal but their Jordan degree types are different. It is stated in the same paper that this
Hilbert function has the smallest socle degree in codimension three, for which there is a pair of AG algebras having the same property. However, in codimension greater than three there are pairs of such AG examples of smaller socle degree \cite[Example~4.4]{Al3}. \par
However, in codimension two, the Jordan degree type of an Artinian graded algebra is determined by the Jordan type, a result that follows from J. Brian\c{c}on's vertical strata in \cite{Bri}.
\end{remark}
\begin{lemma}[Codimension two JDT]\label{cod2JDTlem}\cite[Lemma 2.30]{IMM2} When $A$ is a standard graded algebra of codimension two, and $A$ has Jordan type ${P_{\ell,A }= (p_1, p_2,\ldots,p_s)}$, ${p_1\ge \cdots \ge p_s}$ with respect to an element ${\ell\in A_1}$ and ${\cha \k = 0}$ or ${\cha \k > j}$, the socle degree of $A$, then the Jordan degree type satisfies
$$ \mathcal S_{\ell,A}=\left( (p_1,0),(p_2,1),\ldots , (p_k,k-1),\ldots, (p_s, s-1)\right).$$ 
\end{lemma}
\subsection{Avatars of Jordan degree type.}\label{avatarsec}
The information in the Jordan degree type $\mathcal S_{\ell,A}$ is the same as the numerical information in the central simple modules with respect to $\ell$ of T. Harima and J. Watanabe \cite{HW1,HW2}, (see \cite[Lemma 2.34]{IMM2}). We now show this is equivalent to the information in the ranks of the maps by powers of $\ell$ from the graded pieces of $A$. Recall from Equation \eqref{listJDTeq} that $\eta(p,\nu)$ is the multiplicity of the string $(p,\nu)$ in the JDT $\mathcal S$.\par
\begin{proposition}
\label{JDTrankprop} 
Let $A$ be a standard graded Artinian algebra of socle degree $j$. Then the Jordan degree type $\mathcal S_{\ell,A}$ is equivalent to knowing the ranks of each of the multiplication maps $\ell^{u-k}: A_k\to A_u$, for $0\le k<u\le j$.\par
In particular, if $A$ is graded Artinian Gorenstein the JDT $ \mathcal S_{\ell,A}$ is determined by the set of ranks of all the mixed Hessians, and vice-versa.
\end{proposition}
\begin{proof} 
The second statement concerning mixed Hessians follows from the first and Theorem~\ref{rankthm}. Now, given the JDT, we determine the ranks of maps $\ell^{u-k}: A_k\to A_u$ first by order of $k$, then $u$. For $k=0$ the rank of $\ell^u: A_0\to A_u$ is just the number of strings beginning in degree $0$ and ending in degree at least $u$. For $k>0$, we have
\[
\rk(\ell^{u-k}: A_k\to A_u )=\#\{(p,\nu)\in \mathcal S_{\ell,A}\mid \nu\le k, p> u-\nu\}.
\]
For the converse, it is easy to see that for $0\le \nu\le j$ and $1\le p\le j+1-\nu$ the component $(p,\nu)$ appears $\eta(p,\nu)$ many times in $\mathcal{S}_{\ell,A}$ where
\begin{align*}
\eta(p,\nu)&=\rk(\ell^{p-1}:A_\nu\rightarrow A_{p+\nu-1})\\
&\quad-\big[\rk(\ell^{p}:A_{\nu-1}\rightarrow A_{p+\nu-1})+\rk(\ell^{p}:A_\nu\rightarrow A_{p+\nu})-\rk( \ell^{p+1}:A_{\nu-1}\rightarrow A_{p+\nu})\big].
\end{align*}
\end{proof}
\par
For AG algebras the JDT is symmetric - this was first shown by T. Harima and J. Watanabe  in the context of properties of central simple modules \cite{HW1}, then noted by B. Costa and R. Gondim  \cite[Lemma 4.6]{CGo}; for a proof see also \cite[Proposition 2.38]{IMM2}. This property is useful in determining the possible JDT given a Hilbert function, and is a strong reason for using JDT instead of just JT for graded AG algebras in codimension at least three.\cite{A-N,AAIY}.\par\vskip 0.2cm
\begin{proposition}\label{JDTsymprop}(Symmetry of JDT) Let $A$ be a graded AG algebra and $\ell\in A_1$, and write in list manner 
$\mathcal S_{\ell,A}=\sum_{p}\sum_\nu (p,\nu)^{\eta(p,\nu)}$ (Equation \eqref{listJDTeq}). Then we have for $\nu\le j/2$
\begin{equation}
\eta(p,\nu)=\eta(p,j+1-\nu-p).
\end{equation}
\end{proposition}
\begin{remark}\label{dualrem}(Using dual strings to calculate JDT)
The JDT, like the JT for a graded Artinian algebra may be calculated by considering dual strings in $A^\vee=\Hom(A,\k)=R\circ F$. 
Here  $R=\k[x,y,z]$.\par\vskip 0.2cm
If $A$ is standard graded, and $\ell$ nilpotent, then $A^\vee$ has a decomposition as a direct sum 
$$A^\vee=\sum_i  k[\ell]\circ v_i \text { of simple  $k[\ell]$ modules}.$$
We set $p_i=\dim_\k  k[\ell]\circ v_i$ and recall that the Jordan degree type $\mathcal S(A)=\{(p_i,\deg v_i)^{\eta(p_i,\deg v_i)}\}$, where the pairs   $(p_i,\deg v_i)$ are distinct and $\eta(p_i,\deg v_i)$ is the multiplicity.\par
 Let $A=R/I$ be a standard graded Artinian $\k [x]$-module, $R=\k[x_1,\ldots ,x_r]$ and let the action of $x$ have degree one. Denote by $\overrightarrow{s}: (p,v)$ a string-simple $k[x]$-module $(p,v)$. If $\overrightarrow{s}\subset A$, we may write it as $(p,v)=\{v,xv,\ldots, x^{p-1}v\}$ in degrees $(\deg v,\ldots,\deg v+p-1)$; we let $\overrightarrow{ s}^\vee$ denote an abstract dual string  $(w, x\circ w,\ldots, x^{p-1}\circ w)$, in degrees $\deg v+p-1,\ldots, \deg v$ in the $\k[x]$ module $\k_{DP}[X_1,\ldots, X_r]$, upon which $x$ acts as a contraction lowering degree.  
\begin{lemma}\label{reflectlem}(Using $\ell$-strings in $A^\vee$ to find JDT)  Using the notation above 
let the graded Artinian $\k[x]$-module $A$ satisfy the isomorphism  $A\cong \oplus_i (\overrightarrow {s}_i)$ where $\overrightarrow {s}_i$ are the simple graded dual $k[x]$-modules.  Then we have that $\Hom(A,k)\cong \oplus_i  {\overrightarrow {s}_i}^\vee$, a direct sum of graded simple $\k[x]$-modules (strings of the same length).  Assume  $\mathcal S(A)=\oplus_i \{(p_i,n_i)^{\eta(p_i,n_i)}\}$. Then we have
\begin{equation}
\label{reflecteq} 
\mathcal{S}^\vee(A)=\{(p_i,n_i+(p_1-1))^{\eta(p_i,n_i)}\}.
\end{equation}
\end{lemma}
\begin{proof} 
Writing $A$ as a finite direct sum of simple $\k[x]$-modules, we have $\Hom(A,\k)$ is the analogous finite direct sum of their duals.
\end{proof}
\begin{example} 
The Artinian algebra $A=\k[x,y]/(x^3,y^3,x^2y^2)$ has Hilbert function $H=(1,2,3,2,0)$ and $x$-strings
$\{1,x,x^2\}$, $\{y,yx,yx^2\}$, and $\{y^2,y^2x\}$ so $A$ has SJT $(3_0,3_1,2_2)$.  Here $A=R/I$, $I=\Ann(X^2Y,XY^2)$, whose $x$-strings are  $\{X^2Y,XY,Y\}$, $\{X^2,X,1\}$, and $\{XY^2,Y^2\}$ so $\mathcal {S}^\vee(A)=(3_3,3_2,2_3)$ which for socle degree $3$ corresponds to $\mathcal {S}(A)$, by Equation \eqref{reflecteq}.
\end{example}
\end{remark}
\begin{example}
\cite{AAIY} Let $j\ge 4$, take $F=X^j+X^2Y^{j-2}+XZY^{j-2}$, let $A=R/\Ann F$. Then
 $\Ann F=I=(z^2,x^2z,x^2y-xyz,y^{j-1}, x^{j-1}-y^{j-2}z)$ and  the Hilbert function $H(A)=(1,3,5^{j-3},3,1). $\footnote{To verify the Hilbert function $H(A)$ it helps that according to \cite{BuEi,St,Di} (see \cite[Thm. 5.25]{IK}) that a codimension three Gorenstein ideal of order $\kappa$ can have at most $2\kappa+1$, generators, the given ideal is in $\Ann F$, and the first three generators of $I$ define a length 5 scheme, of Hilbert function $(1,3,\overline{5})$.\par
 Warning: Dual strings are not usually obtained by just changing variables from $x,y,\ldots$ to $X,Y,\ldots$!\par }
 We take $\ell=x$ and claim that $\mathcal{S}_{x,A}=\bigl((j+1)_0,(3\uparrow_1^{j-3}),2_{1,j-2},(1\uparrow_2^{j-2})\bigr)$.\par
 \noindent There are corresponding strings in $A^\vee$ with (dual) generators as $\k[x]$-modules, $(j+1)_0$ corresponds to $F$. For every $1\le i\le j-3$ we have $y^i\circ F = X^2Y^{j-i-2}+XZY^{j-i-2}$ which explains $(3\uparrow_1^{j-3})$. Now let 
$(y^{j-2}-x^{j-2})\circ F=XZ$ and $z\circ F=XY^{j-2}$ so there are strings of length two in degrees $1$ and $j-2$; i.e. $2_{1,j-2}$ exists in $\mathcal{S}_{x,A}$. That leaves $(1\uparrow_2^{j-2})$ to be explained. The difference $H(A)$ minus the HF of the Jordan type so far is just $(0,0,1,1, \ldots, 1_{j-2})$. Since there can be no further strings in $A^\vee$ of length $2$ (which require an $X\cdot h(Y,Z)$ dual generator), the symmetry of JDT forces that $(1\uparrow_2^{j-2})$ is the last set of $j-3$ strings. 
\end{example}
Note that the limit of a family of graded Gorenstein algebras need not be Gorenstein, and need not
have symmetric Jordan degree type, even for simple Hilbert function as $H=(1,3,3,1)$ (\cite[Example 4.16]{AEIY}.\par\vskip 0.3cm
\begin{question} Can we use Lemma \ref{JDTdeformlem} to show the existence of irreducible components in the constant Jordan type strata of an algebra $A$, or in a family of Artinian algebras? \par
\end{question}

\section{Non-graded algebras, and further invariants.}\label{localsec}
One nice application of Jordan type is to show that certain families of non-graded Artinian Gorenstein algebras with given Hilbert function have several irreducible components (Section \ref{3.1sec}). Let $H$ be a \emph{Gorenstein sequence}, i.e.\ a sequence of integers that is the Hilbert function of some AG algebra (not necessarily graded). Let $\Gor(H)$ be the family of AG algebras whose Hilbert function is $H$. Since Jordan type is semi-continuous, if $A$ and $B$ are two algebras in $\Gor(H)$ whose generic Jordan types satisfy ${P_A>P_B}$ then $A$ cannot lie on the border of an irreducible component of $\Gor(H)$ whose general element has Jordan type $P_B$. On the other hand, the symmetric decomposition of the Hilbert function of an AG algebra provides another semi-continuous invariant, namely $N_{i,b}=\dim_\kk \bigl(\mathfrak{m}^i/\bigl(\mathfrak{m}^i\cap (0:\mathfrak{m}^b)\bigr)\bigr)$ (see \cite[Lemma 4.1A and Lemma 4.2A]{I3}). Combining both invariants, in \cite{IM2} the second and third authors of the present survey obtain results on the reducibility of certain families of AG algebras, including the following two, Theorems \ref{threecompthm} and \ref{infseriesthm} in Section \ref{3.1sec} showing that there are infinite collections of such families. \vskip 0.2cm\par
\subsubsection{Symmetric decomposition}
We recall the symmetric decomposition of the associated graded algebra to an AG algebra, from \cite{I3,IM1}. Given an AG algebra $A$, of socle degree $j$, with associated graded algebra $A^*$ we define an ideal $C_A(a)$ componentwise, by  (here $\rho$ is the projection to $A_i$)
$$C_A(a)_i=   \rho\left(  \m^i_A\cap (0:\m_A^{j-a-i})/(\m_A^{i+1}\cap (0:\m_A^{j+1-a-i}))\right).$$
\begin{lemma}\label{QaLemma} Let $A$ be an AG algebra over a field $\k$,  of socle degree $j$. The $A^\ast$, its associated graded algebra, has a filtration by ideals
$$  A^\ast=C_A(0)\supset C_A(1)\subset \cdots \supset C_A(j-2), \text{ such that } Q_A(a)=C_A(a)/C_A(a+1)$$
is a reflexive $\k$-module satisfying
$$Q_A(a)_i\cong Q_A(a)_{j-a-i}^\vee=\Hom_{\k} (Q(a)_{j-a-i},\k).$$
Let a Macaulay dual generator of $A$ be $F_A=F_j+F_{j-1}+\cdots$. Then $Q(0)$ is a socle degree $j$ graded Artinian algebra, whose Macaulay dual generator is $F_j$, and $Q(a)$ is determined by $(F_A)_{\ge j-a}=F_j+\cdots+F_{j-a}$.\footnote{Warning: when $a\ge 1$ the relation between $Q(a)$ and $(F_A)_{\ge j-a}$ is subtle: see \cite[\S 1.3]{IM1}.}\par
Let $H_A(a)=H(Q_A(a))$.  Then $H_A(a)$ is symmetric about $(j-a)/2$, and the Hilbert function $H(A)$ satisfies 
$$H(A)=\sum_a H_A(a).$$
\end{lemma}
\begin{example}\label{decompex}(i). Let  $R=\k [x,y,z,w]$ and $F=Y^6+X^4Z+W^3$; then the Artinian Gorenstein algebra $A=R/\Ann F$ satisfies 
\[
A=(xy,\,xw,\,yz,\,yw,\,zw,\,z^2,\,w^3-x^4z,\,w^3-x^5,\,w^3-y^6) \text { and }Q(0)=R/(x,z,w,y^7)\cong \k[y]/(y^7),
\]
The Macaulay dual is 
\[
A^\vee=R\circ F=\langle F, \{Y^i, 0\in [0,5]\},\{X^iZ,i\in [0,3]\}, \{X^i, i\in [1,4]\}, W,W^2\rangle,
\] 
and we can identify ${Q(1)^\vee=\langle X,Z,X^2,XZ,X^3,X^2Z,X^4,X^3Z\rangle}$, ${Q(2)=0}$, $Q(3)^\vee=\langle W,W^2\rangle$. The symmetric Hilbert function decomposition $\mathcal D(A)$ of $H(A)=(1,4,4,3,3,1,1)$ is thus
\[
\mathcal D(A)=\bigl(H(0)=(1^6),\,H(1)=(0,2,2,2,2,0),\, H(3)=(0,1,1,0)\bigr)
\]
It is not hard to show that this is the unique symmetric decomposition possible for $H(A)$.\par \vskip 0.15cm\noindent 
(ii). In contrast, there are two symmetric decompositions for $H(C)=(1,3,3,2,1)$ of Example~\ref{generalexample1}. For both $C$ and its deformation $C(t)$ of the Example we  have 
\[
H(0)=(1,2,2,2,1),\, H(1)=(0,1,1,0).
\]
But for $G+(X+Y)^4=X^3Y+(X+Y)^4 + Y^2Z$
we have $H(0)=(1,2,3,2,1)$, $H(1)=0$ and $H(2)=(0,1,0)$.\footnote{Note, ${G(t)+(X+Y)^4}$, ${t\neq 0}$ has larger HF $(1,3,4,2,1)$, with decomposition $(1,2,3,2,1)+(0,1,1,0)$.}
\par 
The generic Jordan type for ${\mathcal{R}/\Ann\bigl(G+(X+Y)^4\bigr)}$ is ${(5,3,1,1)}$, the same as for $G$. \par

Note that in the first dual generator $G+(X+Y)^4+Y^2Z$  in (ii) above $Y^2Z$ is an \emph{exotic summand}, since the multiplier $Y^2$ of the new variable $Z$ is a partial of $G+(X+Y)^4$; but here this does not hide $Z$ enough to change the Hilbert function. Exotic summands are defined and discussed in \cite[\S 2.2]{IM1} and \cite{BJMR}, we will not treat them here.\par
\end{example}
Symmetric decompositions of AG algebras, discovered in 1985 (\cite{I2}), have seen increasing appllications recently, particularly to issues of the scheme or cactus length of forms -a study begun by Alessandra Bernardi and Kristen Ranestad in their much-cited 2013 \cite{BR}, see also \cite{BJMR}. They have been applied also to classification of Gorenstein local algebras, as \cite{JMR}.
\subsection{Families $\Gor(H)$ with multiple irreducible components.}\label{3.1sec}
As mentioned earlier the following two results depend on a comparison of the dominance partial order on Jordan types $P_{(\ell,A)}$ with a contrasting partial order arising from the semicontinuity of certain invariants $N_{i,b}=\dim_\kk \bigl(\mathfrak{m}^i/\bigl(\mathfrak{m}^i\cap (0:\mathfrak{m}^b)\bigr)\bigr)$ of symmetric decompositions. Families $\Gor(H)$ with multiple irreducible components had been previously exhibited using the semicontinuity of symmetric decompositions and other arguments, such as dimension, as \cite[Theorem 4.3]{I3} for $H=(1,3,3,2,1,1)$.
\begin{theorem}
\label{threecompthm}
\cite[Theorem 3.3]{IM2}.
Let ${k\ge 2}$ and consider the Gorenstein sequence \linebreak 
${H(k)=(1,3,4^k,3,2,1)}$. Then $\Gor\bigl(H(k)\bigr)$ has exactly three irreducible components, each corresponding to a symmetric decomposition of the Hilbert function $H(k)$.
\end{theorem}

\begin{theorem}
\label{infseriesthm}
\cite[Theorem 3.6]{IM2}.
Let ${k\ge 1}$ and ${s\ge 2}$, and consider the Gorenstein sequence $H(k,s)=(1,3,4^k,3^s,2,1)$. Then $\Gor\bigl(H(k)\bigr)$ has at least two irreducible components, each corresponding to a symmetric decomposition of the Hilbert function $H(k,s)$.
\end{theorem}

\subsection{Sequential Jordan type, L\"oewy Jordan type, Double sequential Jordan type.}
Jordan type is defined for all Artinian algebras; however Jordan degree type does not have a natural extension to non-graded algebras:  Chris McDaniel showed that it is not even possible to find a Jordan basis for the multiplication by $x$ on the AG algebra ${A=\k[x,y]/I}$, ${I=(x^2-xy^2,\,y^4)}$ that is consistent with the Hilbert function ${H(A)=(1,2,2,2,1)}$, so a Jordan order type for non-graded Artinian algebras generalizing JDT for graded Artinians is not possible (\cite[Example 2.15]{IM2}).
We define several refinements of Jordan type introduced in \cite{IMS}, that have some desirable deformation  properties. 
They each depend upon a filtration of an Artinian algebra $A$ by $\m$-adic quotients $A/\m^i$. or by L\"oewy quotients $A/(0:\m^i)$ or by ideals obtained by combining the two.
\begin{definition}[Sequential, L\"oewy and Double Jordan type]\label{SJTdef}
Let $A$ be an Artinian local algebra of socle degree $j$, let $\m$ be its maximal ideal, and let $\ell\in\m$.   
\begin{enumerate}[(i).]
	\item The Sequential Jordan type (SJT) of $(\ell,A)$ is given by the sequence 
	\[ (P_{,\ell,A/{\m^i}}),\ i\in\{1,\ldots,j\} \]
	of Jordan types of successive quotients of $A$ by powers of the maximal ideal.
	\item The L\"oewy Sequential Jordan type (LJT) of $(\ell,A)$ is given by the sequence 
	\[ (P_{\ell,A/{(0:\m^{j-k})}}),\ k\in\{1,\ldots,j\} \]
	of Jordan types of successive quotients of $A$ by the L\"oewy ideals.
	\item The Double Sequential Jordan type  (DSJT) for a pair $(A,\ \ell)\in\m$, is given by the table whose $(a,i)$ entry is the partition
\[ 
P_{\ell,B_{a,i}},\ \text{where}\  B_{a,i}: = A/\bigl(\m^i\cap (0:\m^{j+1-a-i})\bigr),\, 0\le a\le j,\, 0\le i\le j+1-a  
\]
giving the Jordan type of the quotient of $A$ by intersections of a L\"oewy ideal with a power of the maximal ideal.
\end{enumerate}
\end{definition}\noindent

\begin{example}
\label{generalexample3}
Let ${\mathcal{R}=\kk\{x,y,z\}}$ and ${C=\mathcal{R}/\Ann G}$, where ${G=X^3Y + Y^2Z}$, an AG algebra with socle degree ${j=4}$, as in Example \ref{generalexample1}. Then the double sequential Jordan type of an element ${\ell=ax+by+cz+h}$, with ${h\in\m_C^{\,2}}$, and ${ab\ne0}$ is 
\begin{center}
{\scriptsize
\begin{tabular}{cccccc}
$a$, $i$ & $i=0$ & $1$ & $2$ & $3$ & $4$  \\
$a=0$ & & $(1)$ & $(2,1,1)$ & $(3,2,1,1)$ & $(4,3,1,1)$ \\
$1$ & $(1)$ & $(2,1)$ & $(3,2,1)$ & $(4,3,1,1)$ \\
$2$ & $(2,1)$ & $(3,2,1)$ & $(4,3,1,1)$ \\
$3$ & $(3,2,1)$ & $(4,3,1,1)$\\
$4$ & $(4,3,1,1)$  & & & \multicolumn{2}{r}{$P_{\ell,C}=(5,3,1,1)$}
\end{tabular}
}
\\[1em]
DSJT: $P_{\ell,C_{a,i}},\,  C_{a,i}\coloneqq C/\bigl(\m_C^{\,i}\cap (0:\m_C^{\,j+1-a-i})\bigr)$
\end{center}
Note that for ${a=0}$, since ${\m_C^{\,i}\subseteq (0:\m_C^{\,j+1-i})}$, the first row in this table gives us the sequential Jordan type, and for ${i=0}$, since ${\m_C^{\,0}=C}$, the first column shows us the Loewy sequential Jordan type. Also, if ${a+i=j+1}$, we have ${\m_C^{\,i}\cap (0:\m_C^{\,j+1-a-i})=0}$, so the Jordan type is that of the pair $(\ell,C)$. To save space, we write this Jordan type at the lower right\-{-hand} corner of the table. Here, also the diagonal ${a+i=j}$ is constant. This is always the case for an AG algebra, because ${(0:\m_C)=\m_C^{\,j}}$.
\end{example}

\begin{lemma}\cite[Prop. 2.2]{IMS}\label{SJT=JDTlem}
When $A$ is standard graded, both the Sequential Jordan type and the L\"oewy Sequential Jordan type are equivalent to the Jordan degree type.
\end{lemma}
\begin{lemma}\cite[Theorem 2]{IMS} We have that $JT\ge SJT\ge DSJT$ and $JT\ge LSJT\ge DSJT$ are true refinements.
\end{lemma}
We do not know if $SJT+LSJT$ to $LSJT$ is a true refinement. 
Given their definitions, and the dominance order (Definition \ref{dominancedef}) for Jordan type, we can define natural dominance orders for SJT, LSJT, and DSJT when the total lengths of the algebras compared are the same.  For example
an SJT $S$ dominates $S^\prime$ if for every degree $i$, the partition $S_{\le i}\ge S^\prime_{\le i}$ in the dominance order of Definition \ref{dominancedef}. Recall that we take 
$\mathcal R=\k[x_1,\ldots, x_r]$ with maximal ideal $\m=(x_1,\ldots, x_r)$. 
We need a preparatory result, adapted from \cite[Corollary 2.44]{IMM2}, that will imply the deformation results we wish to show for these invariants. \vskip 0.2cm\par

\begin{lemma}[Semicontinuity of Jordan type]\label{semicontlem} 
(i). Let $M(\tau)$ for $\tau\in \Z$ be a family of constant length $R$-modules over a parameter space $\Z$ and let $\ell\in \m_R$. Then for a neighborhood $U_0\subset \Z$ of $\tau_0$, we have that $\tau\in U_0\Rightarrow P_{\ell,M(\tau)}\ge P_{\ell, M(\tau_0)}$ in the dominance partial order of Definition~\ref{dominancedef}.\par
(ii). Let $M$ be a finite module over an Artinian algebra $A$ and $\{\ell(t), t\in \mathfrak X\}$, $\mathfrak X$ a curve, be a family of linear forms or of elements of $\m$ (according to whether $A$ is graded or local). 
Assume that $P_{\ell,M}=P$ is constant for $\ell\in \mathcal U\subset \mathfrak X$, an open dense in $\mathfrak X$; and let $\ell_0\in \mathfrak X\backslash U$. Then $P\ge P_{\ell_0,M}$ in the dominance partial order.
\end{lemma}

Since an Artinian algebra is a deformation of its associated graded algebra, we have
\begin{corollary}\label{AGAcor}  for any Artinian $A$, and any $\ell\in \m_A$, the Jordan type $P_{\ell,A}$ of $A$ dominates the Jordan type $P_{\ell,A^\ast}$ of the associated graded algebra of $A$.
\end{corollary} 
See Examples \ref{1ex} and \ref{generalexample1}.
\begin{proposition}[Deformations]\label{SJTdeformprop}\cite[Prop. 2.12]{IMS}
Let $A(\tau), \tau\in \Z$ be a flat (constant length) family of Artinian algebra quotients of $R$, let $\tau_0\in \Z$, and fix $\ell\in \mathcal R$. 
\begin{enumerate}[(i).]
\item (SJT) Assume that the Hilbert function $H(A(\tau))$ is constant. Then there is an open neighborhood $\mathcal U$ of $\tau_0$ in $\Z$ such that $\tau\in \mathcal U$ implies that the SJT $\mathcal S(\ell,\tau)$ of the pair $(\ell, A(\tau))$ dominates the SJT $\mathcal S(\tau_0)$ of $(\ell, A(\tau_0))$
\item (LSJT) Assume that the dimensions of the L\"oewy ideals ${(0:\m_{A(\tau)}^{\,i})}$ are constant along the family $\Z$. Then there is an open neighborhood $\mathcal U$ of $\tau_0$ such that $\tau\in \mathcal U$ implies that the LSJT for the pair $(\ell,A(\tau))$ dominates the LSJT for $(\ell, A(\tau_0)$.
\item (DSJT) Assume that for each pair $(i,k)$ the dimensions of the ideals ${\m_{A(\tau)}
^{\,i}\cap(0:\m_{A(\tau)}^{\,k})}$ are constant along the family $\Z$. 
Then there is an open neighborhood $\mathcal U$ of $\tau_0$ such that $\tau\in \mathcal U$ implies that the DSJT for the pair $(\ell,A(\tau))$ dominates the DSJT for $(\ell, A(\tau_0)$.
\end{enumerate}
\end{proposition}
\begin{proof} The proof is immediate from the semicontinuity of the appropriate invariant on the given locus.
\end{proof}\par
\begin{question} There are many other ways for an Artinian Gorenstein to use the symmetric decomposition to make up a Jordan type related invariant. For example, we could consider the Jordan degree type of each component $Q_A(a)$ of the symmetric decomposition of $A^\ast$, ${a\in [0,j-2]}$, where $j$ is the socle degree. These JDT will satisfy appropriate deformation properties, for example if one fixes the symmetric decomposition of the Hilbert function. Which of these can we use to study irreducible components of symmetric decomposition strata $\mathcal D=\bigl( H(Q(0)), H(Q(1)), \ldots,  H(Q(j-2))\bigr)$ or of $\Gor(H)$, the family of Gorenstein local algebras of a given Hilbert function?
\end{question}
\section{Open questions.}
\begin{question} Can we use SJT, LJT, DSJT to help show that certain families of local Artinian algebras of, say, fixed Hilbert function, have several irreducible components?
\end{question}
\begin{question}
 What is the relation between JT and Betti diagrams? See recent articles of N. Abdallah, and N. Abdallah and H. Schenck \cite{Ab,AbSc}.\par
 \end{question}
 \begin{question}\label{CellsZHquest} Consider Brian\c{c}on's vertical cells (Jordan cells) for the family of codimension two local Artinian algebra of Hilbert function $H$. Answer analogous questions to those of \cite{AIK} involving numbers of generators, and also concerning the symmetric decomposition.
\end{question}
\begin{question}
The following is related to Section \ref{majorsec} Problem viii. The Jordan block decomposition for a (similarity class) of matrices - the Jordan normal form (JNF) has been long known. But only relatively recently - 2008 - had there been work on which pairs of partitions can occur of JNF of two commuting $n\times n$ matrices.\footnote{Early researchers in this area
 seemed more interested in the maximum dimension of a vector space of commuting matrices.} This problem reduces to 
 considering two nilpotent matrices, where JNF is just a partition of $n$. A partition is called
\emph{stable} if its parts differ pairwise by at least two.  One perhaps surprising fact, shown by T. Ko\v{s}ir and P. Oblak with help from 
others is that a commuting pair cannot consist of two different stable partitions: the result has to do with the Hilbert functions of complete intersection quotients of $k[x,y]$.  The problem of finding the maximum (in dominance order) Jordan type (it is stable) commuting with a given partition has attracted some interest, in particular a conjecture of P. Oblak answering this has recently proved by R. Basili.
The problems in this area seem quite difficult, they include study of certain graph associated to a partition; it is as if there is
 some hidden structure lurking behind what we know. See \cite{JS}, and the recent survey by Leila Khatami \cite{Kh} and the references cited there - we have not included references here on this rich topic.
\end{question}

\end{document}